\documentclass[reqno,12pt]{amsart}
\usepackage{amsmath,amssymb}
\usepackage{enumitem}
\usepackage[margin=1.5in]{geometry}

\newtheorem{theorem}{Theorem}
\newtheorem{remark}{Remark}
\newtheorem{example}{Example}

\newcommand{\R}{\mathbb R}
\newcommand{\N}{\mathbb N}
\newcommand{\T}{\mathcal T}
\newcommand{\D}{\mathcal D}
\newcommand{\sT}{\mathfrak T}
\newcommand{\sD}{\mathfrak D}
\newcommand{\sF}{\mathfrak F}
\newcommand{\osp}{\mathfrak{osp}(1|2)}
\newcommand{\Db}{\overline D}

\begin{document}

\title{$\mathfrak{sl}(2)$- and $\mathfrak{osp}(1|2)$-invariant multidifferential operators}

\author[Mtaouaa]{Wafa Mtaouaa}
\address{Department of Mathematics, Faculty of Sciences of Tunis, University of Tunis El Manar, LR Analyse Non-Lin\'eaire et G\'eom\'etrie, LR21ES08, El Manar 2, 2092, Tunis, Tunisia.}
\email{mtawaa.wafa@yahoo.fr}

\date{\today}

\begin{abstract}
We determine the dimension of the space of $\mathfrak{sl}(2)$-invariant $n$-ary differential operators on weighted densities for arbitrary weights, including the singular cases in which the usual elimination procedure breaks down. The proof relies on the Lefschetz property of monomial complete intersections, due to Stanley \cite{s} and Watanabe \cite{w}, applied block by block to the system. In the classical case this recovers, by a different proof and with a single formula valid for all weights, the dimension obtained by Arnal and Ben Ammar \cite{aba}. We also investigate the corresponding problem in the super case.
\end{abstract}

\maketitle

\section{Introduction}

Invariant differential operators associated with Lie algebras and Lie superalgebras play an important role in representation theory and differential geometry and have attracted considerable attention.

 In the framework of weighted densities, a transvectant is an $\mathfrak{sl}(2)$-invariant differential operator acting on tensor products of spaces of densities $\mathcal{F}_{\lambda}$. In the bilinear case, these operators are known as Rankin–Cohen brackets. Their generalization to an arbitrary number of arguments, referred to as $n$-transvectants, was studied by Basdouri, Bartouli, Lerbet and Nasri \cite{bbln}, while the ternary case was investigated by Basdouri and Nasri \cite{be}.

A similar problem arises for the orthosymplectic Lie superalgebra $\mathfrak{osp}(1|2)$ acting on weighted densities over the superspace $\mathbb{R}^{1|1}$. The corresponding density modules and the action of the contact algebra are discussed in \cite{ab}. Invariant bilinear differential operators of this type are known as supertransvectants and were introduced by Gargoubi and Ovsienko \cite{gov}. In this paper, we consider their $n$-ary generalizations, which we call $n$-supertransvectants.

\medskip
 Throughout, $k=\mu-|\lambda|$ with $|\lambda|=\lambda_1+\dots+\lambda_n$, $t_i=-2\lambda_i$, and $\delta_m(t)$ is the number of $\alpha\in\mathbb N^n$ with $\alpha_i\le t_i$ and $|\alpha|=m$. The classical case needs $\delta_m(t)$ and the super case needs $\delta_m(2t)$. We write $\mathbb N=\{0,1,2,\dots\}$.

In the classical setting, the determination of $n$-transvectants amounts to solving a homogeneous linear system, given by \eqref{sys}, whose coefficients depend on the weights $\lambda_1,\ldots,\lambda_n$. The methods developed in \cite{be, bbln} provide a solution when at least one coefficient of the system does not vanish for any relevant value of the index. The work of Basdouri and his coauthors stops at the level of the system \eqref{sys}: they did not give the dimensions of the spaces of invariant operators, and they did not treat the singular cases. However, these methods do not determine the dimension of the solution space in the singular cases, where the usual elimination procedure breaks down. Such cases are not exceptional: they arise whenever all the numbers $-2\lambda_i$ belong to $\{0,\ldots,k-1\}$, where $k=\mu-|\lambda|$. In these situations, the dimension of the space of invariant operators may exceed its generic value.

\medskip
A complete answer to the classical problem is contained, in a different framework, in the work of Arnal and Ben Ammar \cite{aba}. They compute the cohomology spaces $H^j(\mathfrak{sl}(2),\mathcal{D}_{\lambda,\mu})$, $j=0,1,2$, where $\mathcal{D}_{\lambda,\mu}$ denotes the space of multidifferential operators; the invariant operators correspond to $j=0$. These spaces are characterized by a linear map $F:K_k\to K_{k-1}$ between explicit finite-dimensional spaces. When $t=-2\lambda\in\mathbb N^n$ and $k>t_i$ for all $i$, they restrict $F$ to a map $f$ between the spaces spanned by the multi-indices $\alpha$ with $\alpha_i\le t_i$, and prove by induction on $n$ that $f$ has maximal rank. This gives $\dim H^0=\binom{n+k-2}{k}+\delta_{k-1}(t)-\delta_k(t)$ for $\frac{|t|+1}2\le k\le|t|+1$ and the generic value $\binom{n+k-2}{k}$ otherwise \cite[Cor.~11]{aba} (their $d_{nkt}$ is our $\delta_k(t)$). 

\medskip
Our dimension formula for $j=0$ is equivalent to theirs.  Instead of restricting to the sub-basis $\alpha_i\le t_i$, changing basis, passing to the adjoint to reduce to $k\ge\frac{|t|+1}2$, ordering $t_1\le\dots\le t_n$ and inducting on $n$ through a case analysis of block matrices, we split the system \eqref{sys} once and for all into $2^n$ independent blocks, each of which is multiplication by $y_1+\dots+y_n$ on a monomial complete intersection. The Lefschetz property then applies directly to every block, with no induction. Only the block $S=\emptyset$ can contribute beyond the generic count, and it is the transpose of the map $f$ of \cite{aba}; this explains why the part of the space outside the sub-basis can be discarded there.

\medskip
 A single formula, $\dim\T^n_{\lambda,\mu}=\binom{n+k-2}{k}+\max(0,\delta_{k-1}-\delta_k)$, holds for every $k\ge1$ and every $t\in\mathbb N^n$, including $t_i\ge k$, and the range where the correction is non-zero follows from the symmetry and unimodality of $(\delta_m)$ (Remark~1).
 The block decomposition shows which blocks carry the solutions of \eqref{sys} (Remark~3), and not only how many solutions there are.

\medskip
The main results of this paper are divided into two parts.

In the classical case, considered in Sections 2 and 3, we reformulate the system \eqref{sys} as the transpose of an explicit linear map $\psi$. We show that this map decomposes into independent blocks, each of which can be identified with multiplication by $y_1+\cdots+y_n$ on a suitable truncated polynomial algebra. The Lefschetz property of monomial complete intersections, due to Stanley \cite{s} and Watanabe \cite{w} (see also \cite{hmnw}), then provides the rank of each block (the link between the maximal rank of the map $f$ of \cite{aba} and the Lefschetz property was first pointed out in \cite[Remark~1, Sect.~5]{aba}). This leads to an explicit formula for the dimension of the space of $\mathfrak{sl}(2)$-invariant operators, consisting of a binomial coefficient together with a correction term involving the numbers of bounded compositions of certain integers (Theorem 1). In addition to determining the dimension, the block decomposition provides information about the invariant operators themselves.

The second part, developed in Sections 4 and 5, concerns the superalgebra $\mathfrak{osp}(1|2)$. In this setting, the operators $\Omega^\alpha$ used in the classical case are no longer suitable, since they do not distinguish the odd derivatives of the individual arguments. We therefore introduce a different basis, denoted by $D^\alpha$. Using the fact that $\mathfrak{osp}(1|2)$ is generated by its odd part, we reduce the invariance conditions to two equations associated with $X_\theta$ and $X_{x\theta}$. After eliminating the odd coefficients, the problem becomes the computation of the kernel of a single odd operator $B$ acting on
$$
\mathbb{R}[z_1,\ldots,z_n]\otimes
\Lambda(\xi_1,\ldots,\xi_n).
$$
We then decompose this space into suitable blocks and use the representation theory of $\mathfrak{osp}(1|2)$ to determine the dimension of the kernel. The resulting formula, stated in Theorem 2, has a structure analogous to that obtained in the classical case: a generic binomial term together with a correction term arising in the singular cases.

In both settings, the dimension formulas are illustrated through explicit examples and checked by direct computations in low-dimensional cases. The correction terms are non-zero only for a finite range of values determined by the weights.

\section{The classical setting}\label{sec-classical}

Let $\mathcal F_\lambda=\{f\,dx^\lambda : f\in C^\infty(\mathbb R)\}$. A vector field $X_h=h\,\frac d{dx}$ acts on it by $L^\lambda_{X_h}(f\,dx^\lambda)=(hf'+\lambda h'f)\,dx^\lambda$. For $\lambda=(\lambda_1,\dots,\lambda_n)\in\R^n$ and $\mu\in\R$, let $\D_{\lambda,\mu}$ be the space of $n$-ary differential operators $\mathcal F_{\lambda_1}\otimes\dots\otimes\mathcal F_{\lambda_n}\to\mathcal F_\mu$, on which $X_h$ acts by $X_h\cdot A=L^\mu_{X_h}\circ A-A\circ L^\lambda_{X_h}$, with $L^\lambda_{X_h}$ extended by the Leibniz rule. The algebra $\mathfrak{sl}(2)=\mathrm{Span}(X_1,X_x,X_{x^2})$ acts by restriction, and we write
$$
\T^n_{\lambda,\mu}=\D_{\lambda,\mu}^{\mathfrak{sl}(2)}.
$$

Every operator is a sum $\sum_\alpha A_\alpha(x)\Omega^\alpha$, $\alpha\in\N^n$, where
$$
\Omega^\alpha(f_1dx^{\lambda_1}\otimes\dots\otimes f_ndx^{\lambda_n})=\frac1{\alpha!}f_1^{(\alpha_1)}\cdots f_n^{(\alpha_n)}dx^\mu .
$$
Put $k=\mu-|\lambda|$. A direct computation gives
$$
X_h\cdot \sum A_\alpha\Omega^\alpha=\sum_\alpha\Big(A_\alpha' h+(k-|\alpha|)A_\alpha h'\Big)\Omega^\alpha-\frac12\sum_\alpha\sum_{\alpha_i>0}(\alpha_i-1+2\lambda_i)A_\alpha h''\,\Omega^{\alpha-\varepsilon_i}.
$$
Since $X_x=\frac12[X_1,X_{x^2}]$, an operator is invariant if and only if it is killed by $X_1$ and $X_{x^2}$. Being killed by $X_1$ means that the $A_\alpha$ are constants. Then $X_x\cdot\Omega^\alpha=(k-|\alpha|)\Omega^\alpha$, so $\T^n_{\lambda,\mu}=0$ unless $k\in\N$, and for $k\in\N$ an invariant operator is a combination of the $\Omega^\alpha$ with $|\alpha|=k$. For $k=0$ the only such operator is $\Omega^{0}(f_1dx^{\lambda_1}\otimes\dots\otimes f_ndx^{\lambda_n})=f_1\cdots f_n\,dx^\mu$; it is invariant (the system below is empty), so $\dim\T^n_{\lambda,\mu}=1$. Assume $k\ge1$ from now on. Writing $t_i=-2\lambda_i$, the condition $X_{x^2}\cdot A=0$ for $A=\sum_{|\alpha|=k}a_\alpha\Omega^\alpha$ becomes
\begin{equation}\label{sys}
\sum_{i=1}^n(\alpha_i-t_i)\,a_{\alpha+\varepsilon_i}=0\qquad(|\alpha|=k-1).
\end{equation}
This is a system of $\binom{n+k-2}{k-1}$ equations in $\binom{n+k-1}{k}$ unknowns, and $\dim\T^n_{\lambda,\mu}$ is the number of unknowns minus the rank of the system.

\section{The classical result}\label{sec-classical-result}

For $t\in\N^n$ and $m\in\mathbb Z$ put
$$
\delta_m(t)=\#\{\alpha\in\N^n:\ \alpha_i\le t_i\ (1\le i\le n),\ |\alpha|=m\},
$$
so that $\sum_m\delta_m(t)q^m=\prod_i(1+q+\dots+q^{t_i})$, and $\delta_m(t)=0$ for $m<0$. In this section $t=-2\lambda$ and we write $\delta_m=\delta_m(t)$.

\begin{theorem}\label{main}
Let $k=\mu-|\lambda|$. Then $\T^n_{\lambda,\mu}=0$ unless $k\in\N$, and for $k=0$ one has $\dim\T^n_{\lambda,\mu}=1$, the space being spanned by $\Omega^0$ (this is again $\binom{n+k-2}{k}$). Assume now $k\in\N$, $k\ge1$.
\begin{enumerate}[label=(\roman*)]
\item If $-2\lambda_i\notin\{0,\dots,k-1\}$ for some $i$, then $\dim\T^n_{\lambda,\mu}=\binom{n+k-2}{k}$. If for instance $i=1$, the map $A=\sum a_\alpha\Omega^\alpha\mapsto(a_\alpha)_{\alpha_1=0}$ is an isomorphism from $\T^n_{\lambda,\mu}$ onto $\R^{\binom{n+k-2}{k}}$, and a basis is $\{A_\gamma:\ \gamma_1=0,\ |\gamma|=k\}$, where $A_\gamma$ is the invariant operator with $a_\alpha=\delta_{\alpha\gamma}$ for all $\alpha$ with $\alpha_1=0$. (In general $\Omega^\gamma$ itself is not invariant.)
\item If $t=-2\lambda\in\{0,\dots,k-1\}^n$, then
$$
\dim\T^n_{\lambda,\mu}=\binom{n+k-2}{k}+\max\big(0,\ \delta_{k-1}-\delta_k\big).
$$
\end{enumerate}
\end{theorem}

\begin{proof}
(i) Say $-2\lambda_1\notin\{0,\dots,k-1\}$. Then $\alpha_1-t_1\ne0$ for every $\alpha_1\le k-1$. Choose the $a_\alpha$ with $\alpha_1=0$ freely. For $\beta$ with $\beta_1=0$ and $|\beta|=k-1$, the equation \eqref{sys} for $\alpha=\beta$ expresses $a_{\beta+\varepsilon_1}$ through the $a_{\beta+\varepsilon_i}$, $i\ge2$, which have first index $0$. Next, for $j=1,\dots,k-1$ the equations \eqref{sys} with $\alpha_1=j$ determine the $a_\alpha$ with $\alpha_1=j+1$ from those with $\alpha_1=j$ (the coefficient $j-t_1$ is non-zero); these equations and these unknowns correspond bijectively through $\beta\mapsto\beta+\varepsilon_1$, so no further condition appears. Hence the free parameters are exactly the $a_\alpha$ with $\alpha_1=0$, there are $\binom{n+k-2}{k}$ of them, and $A\mapsto(a_\alpha)_{\alpha_1=0}$ is an isomorphism.

(ii) We compute the rank of \eqref{sys} in four steps.

\smallskip
\emph{Step 1: transposition.} Let $V_m$ be a vector space with basis $(e_\beta)_{|\beta|=m}$ and define $\psi:V_{k-1}\to V_k$ by $\psi(e_\beta)=\sum_i(\beta_i-t_i)e_{\beta+\varepsilon_i}$. The matrix of \eqref{sys} is the transpose of that of $\psi$, so the rank of the system is $\operatorname{rank}\psi$. It is convenient to let $\psi=E_1+\dots+E_n$ act on $V=\bigoplus_mV_m$, where $E_i(e_\beta)=(\beta_i-t_i)e_{\beta+\varepsilon_i}$ only moves the $i$-th index.

\smallskip
\emph{Step 2: cutting the chains.} In one variable, $E(e_j)=(j-t)e_{j+1}$ has a zero coefficient exactly at $j=t$, so the chain $e_0\to e_1\to\cdots$ is cut between $e_t$ and $e_{t+1}$. Rescale: let $\kappa_j=\prod_{r<j}(r-t)$ for $j\le t$, $\kappa_j=(j-t-1)!$ for $j\ge t+1$ (all non-zero), and $u_j=\kappa_je_j$. Then
$$
E(u_j)=u_{j+1}\ (j\ne t),\qquad E(u_t)=0.
$$
Doing this in each variable (with $t=t_i$) we get a basis $u_\beta=\kappa_\beta e_\beta$ of $V$ in which each $E_i$ sends $u_\beta$ to $u_{\beta+\varepsilon_i}$, except when $\beta_i=t_i$, where it gives $0$.

\smallskip
\emph{Step 3: blocks.} For $S\subset\{1,\dots,n\}$ let
$$
\Gamma_S=\{\beta\in\N^n:\ \beta_i\ge t_i+1\ (i\in S),\ \beta_i\le t_i\ (i\notin S)\},
$$
let $V_{S,m}$ be spanned by the $u_\beta$ with $\beta\in\Gamma_S$, $|\beta|=m$, and $g_S(m)=\dim V_{S,m}$. The sets $\Gamma_S$ partition $\N^n$, so $V_m=\bigoplus_SV_{S,m}$ and $\sum_Sg_S(m)=\binom{n+m-1}{m}$. By Step~2, $\psi$ maps $V_{S,k-1}$ into $V_{S,k}$: for $\beta\in\Gamma_S$, $E_i(u_\beta)$ is either $0$ or $u_{\beta+\varepsilon_i}$ with $\beta+\varepsilon_i\in\Gamma_S$. Moreover
$$
\psi(u_\beta)=\sum_{i:\ \beta+\varepsilon_i\in\Gamma_S}u_{\beta+\varepsilon_i}.
$$
Hence $\operatorname{rank}\psi=\sum_S\operatorname{rank}\psi_S$ with $\psi_S=\psi|_{V_{S,k-1}}$.

\smallskip
\emph{Step 4: Lefschetz.} Put $\sigma_S=\sum_{i\in S}(t_i+1)\varepsilon_i$, $d_S=|\sigma_S|$ and
$$
\mathcal B_S=\bigotimes_{i\notin S}\R[y_i]/(y_i^{t_i+1})\ \otimes\bigotimes_{i\in S}\R[y_i]/(y_i^{k+1}),\qquad \ell=y_1+\dots+y_n.
$$
The map $u_\beta\mapsto y^{\beta-\sigma_S}$ is a bijection from $\{u_\beta:\beta\in\Gamma_S,|\beta|=m\}$ onto the monomial basis of the degree $m-d_S$ part of $\mathcal B_S$, for $m\le k$: indeed the exponents are $\ge0$, those with $i\notin S$ are $\le t_i$, and those with $i\in S$ are $\le k$ when $|\beta|\le k$. Multiplication by $y_i$ on $y^{\beta-\sigma_S}$ is non-zero exactly when $\beta+\varepsilon_i\in\Gamma_S$ (for $i\notin S$ this needs $\beta_i<t_i$, for $i\in S$ the exponent stays below $k+1$). So under this identification $\psi_S$ is multiplication by $\ell$ from degree $k-1-d_S$ to degree $k-d_S$ of $\mathcal B_S$. By the Lefschetz property of monomial complete intersections, due to Stanley and Watanabe \cite{s,w} (see also \cite{hmnw}), multiplication by $\ell$ on $\R[y]/(y_1^{a_1},\dots,y_n^{a_n})$ has maximal rank in every degree (only this weak form is used here); this is the fact that \cite[Remark~1, Sect.~5]{aba} uses for their map $f$. Therefore $\operatorname{rank}\psi_S=\min(g_S(k-1),g_S(k))$, and
$$
\dim\T^n_{\lambda,\mu}=\sum_S\big(g_S(k)-\min(g_S(k-1),g_S(k))\big)=\sum_S\max\big(0,g_S(k)-g_S(k-1)\big).
$$

It remains to count. Writing $\beta=\sigma_S+\gamma$ we get $g_S(m)=[q^{m-d_S}]H_S(q)$ with
$$
H_S(q)=\prod_{i\notin S}(1+q+\dots+q^{t_i})\,(1-q)^{-|S|},
$$
and $g_\emptyset=\delta$. For $S\ne\emptyset$ the series $(1-q)H_S(q)=\prod_{i\notin S}(1+\dots+q^{t_i})(1-q)^{-(|S|-1)}$ has non-negative coefficients, so $g_S(k)-g_S(k-1)\ge0$ and the maximum is not needed. Then
$$
\sum_S\big(g_S(k)-g_S(k-1)\big)=\binom{n+k-1}{k}-\binom{n+k-2}{k-1}=\binom{n+k-2}{k},
$$
and the term $S=\emptyset$ of this sum equals $\delta_k-\delta_{k-1}$. Replacing it by $\max(0,\delta_k-\delta_{k-1})$ gives
$$
\dim\T^n_{\lambda,\mu}=\binom{n+k-2}{k}-(\delta_k-\delta_{k-1})+\max(0,\delta_k-\delta_{k-1})=\binom{n+k-2}{k}+\max(0,\delta_{k-1}-\delta_k).\qedhere
$$
\end{proof}

\begin{remark}
(1) The sequence $(\delta_m)$ is symmetric, $\delta_m=\delta_{\tau-m}$ with $\tau=|t|=-2|\lambda|$ (use $\alpha\mapsto t-\alpha$), and unimodal: for $2m<\tau$ the map $\ell^{\tau-2m}:\mathcal A_m\to\mathcal A_{\tau-m}$ on $\mathcal A=\R[y]/(y_i^{t_i+1})$ is bijective by the strong Lefschetz property, hence $\ell:\mathcal A_m\to\mathcal A_{m+1}$ is injective. Consequently the correction in (ii) is positive only when $\frac{\tau+1}2<k\le\tau+1$; for $2k\le\tau+1$ or $k\ge\tau+2$ the dimension is the generic one, $\binom{n+k-2}{k}$.\\
(2) The proof of (ii) works for every $t\in\N^n$. If some $t_i\ge k$ it gives again $\delta_{k-1}\le\delta_k$, hence part (i) for $t\in\N^n$ with some $t_i\ge k$.\\
(3) The system \eqref{sys} splits into one independent system for each block $S$, which gives explicit solutions. For $S\ne\emptyset$ with $d_S=k$ the source of $\ell$ is zero, and $\Omega^{\sigma_S}$ itself is invariant.
\end{remark}

\begin{remark}[Comparison with \cite{aba}]
In the notation of \cite{aba}, $d_{nkt}=\delta_k(t)$ and $h^0_{nkt}=\dim H^0(\mathfrak{sl}(2),\D_{\lambda,\mu})=\dim\T^n_{\lambda,\mu}$. Put $\tau=|t|$. By Remark~1, $\max(0,\delta_{k-1}-\delta_k)$ equals $\delta_{k-1}-\delta_k$ for $\frac{\tau+1}2\le k\le\tau+1$ and $0$ otherwise, so Theorem~\ref{main}\,(ii) is exactly \cite[Cor.~11]{aba}. The block $S=\emptyset$ is, up to the rescaling of the basis in Step~2 and transposition, the map $f$ of \cite{aba}. For $S\neq\emptyset$ the proof shows that $g_S(k)-g_S(k-1)\ge0$, so these blocks only contribute to the generic term $\binom{n+k-2}{k}$; this is the content of \cite[Lemma~7]{aba}, namely $\operatorname{coker}F\simeq\operatorname{coker}f$.
\end{remark}

\begin{example}\label{ex-classic}
Take $n=2$, $\lambda_1=\lambda_2=-\frac12$, $\mu=1$, so $t=(1,1)$ and $k=2$. Here $(\delta_0,\delta_1,\delta_2)=(1,2,1)$, hence
$$
\dim\T^2_{\lambda,\mu}=\binom22+\max(0,\delta_1-\delta_2)=1+1=2.
$$
Directly, \eqref{sys} for $\alpha=(1,0)$ and $\alpha=(0,1)$ reads $-a_{(1,1)}=0$ in both cases, so $a_{(2,0)}$ and $a_{(0,2)}$ are free and $\T^2_{\lambda,\mu}=\mathrm{Span}\{\Omega^{(2,0)},\Omega^{(0,2)}\}$. In terms of blocks, $\Gamma_\emptyset\ni(1,1)$ contributes $0$ because $\ell:\mathcal B_\emptyset[1]\to \mathcal B_\emptyset[2]$ goes from dimension $2$ onto dimension $1$, while $\Gamma_{\{1\}}\ni(2,0)$ and $\Gamma_{\{2\}}\ni(0,2)$ contribute $1$ each.

For $n=3$, $t=(1,1,1)$ we have $(\delta_m)=(1,3,3,1)$ and the theorem gives $\dim\T^3_{\lambda,\mu}=3,6,6,6,7$ for $k=2,\dots,6$. We also checked these values against the rank of \eqref{sys} computed directly.
\end{example}

\section{The super setting}\label{sec-super}

Let $C^\infty(\mathbb R^{1|1})=\{f_0(x)+f_1(x)\theta\}$, with $\theta$ odd. Put $D=\partial_\theta+\theta\partial_x$ and $\Db=\partial_\theta-\theta\partial_x$, so that $D^2=\partial_x$. The contact vector fields are
$$
X_F=F\partial_x+\tfrac12D(F)\Db,\qquad F\in C^\infty(\mathbb R^{1|1}),
$$
and $[X_F,X_G]=X_{\{F,G\}}$ with $\{F,G\}=FG'-F'G+\frac12D(F)\Db(G)$. We use the subalgebra
$$
\osp=\mathrm{Span}(X_1,X_x,X_{x^2},X_{x\theta},X_\theta).
$$
Its even part is $\mathfrak{sl}(2)$, and its odd part is spanned by $X_\theta=\frac12D$ and $X_{x\theta}=\frac12xD$. The odd generators satisfy
$$
[X_\theta,X_\theta]=\tfrac12X_1,\qquad [X_{x\theta},X_\theta]=\tfrac12X_x,\qquad [X_{x\theta},X_{x\theta}]=\tfrac12X_{x^2},
$$
where the brackets are anticommutators. Hence
\begin{equation}\label{gen}
X_1=4X_\theta^2,\qquad X_x=2\,(X_{x\theta}X_\theta+X_\theta X_{x\theta}),\qquad X_{x^2}=4X_{x\theta}^2 .
\end{equation}

The densities of weight $\lambda$ are $\sF_\lambda=\{F\omega^\lambda\}$, with $\omega=dx+\theta d\theta$, and $X_F$ acts by $\mathfrak L^\lambda_{X_F}(G\omega^\lambda)=(X_F(G)+\lambda F'G)\omega^\lambda$. In particular $\mathfrak L^\lambda_{X_\theta}=\frac12D$ and $\mathfrak L^\lambda_{X_{x\theta}}=\frac12xD+\lambda\theta$. Let $\sD_{\lambda,\mu}$ be the space of $n$-ary differential operators $\sF_{\lambda_1}\otimes\dots\otimes\sF_{\lambda_n}\to\sF_\mu$, with
$$
X_F\cdot A=\mathfrak L^\mu_{X_F}\circ A-(-1)^{|A||F|}A\circ\mathfrak L^\lambda_{X_F},
$$
and let $\sT^n_{\lambda,\mu}=\sD_{\lambda,\mu}^{\osp}$. By \eqref{gen},
$$
\sT^n_{\lambda,\mu}=\ker X_\theta\cap\ker X_{x\theta}.
$$

The operators $\Omega^\alpha$ of the classical case do not work here: they cannot see the odd derivative of each argument separately, and the span of the $\Omega^\alpha$ is not stable under $\osp$ for $n\ge2$. We use instead, for $\alpha\in\N^n$,
$$
D^\alpha(F_1\omega^{\lambda_1}\otimes\dots\otimes F_n\omega^{\lambda_n})=D^{\alpha_1}(F_1)\cdots D^{\alpha_n}(F_n)\,\omega^\mu .
$$
The $D^\alpha$ span $\sD_{\lambda,\mu}$ over $C^\infty(\mathbb R^{1|1})$, so every operator is written $A=\sum_\alpha a_\alpha(x,\theta)D^\alpha$. On the submodule of polynomial operators we use the basis
$$
a_{m,\alpha}=x^mD^\alpha,\qquad b_{m,\alpha}=x^m\theta D^\alpha\qquad(m\in\N,\ \alpha\in\N^n).
$$

\section{The super result}\label{sec-super-result}

Keep $k=\mu-|\lambda|$ and $t=-2\lambda$. For $\alpha\in\N^n$ put $s_i(\alpha)=\sum_{j\ge i}\alpha_j$ and
$$
p_i(\alpha)=\begin{cases}\alpha_i/4 & \alpha_i\ \text{even},\\ (\alpha_i-1)/4+\lambda_i & \alpha_i\ \text{odd},\end{cases}
\qquad
\mathcal E(\beta)=\{i:\ \beta_i\ge1,\ s_i(\beta)\ \text{even}\}.
$$

\begin{theorem}\label{mainsuper}
Let $n\ge2$ and $k=\mu-|\lambda|$. Then $\sT^n_{\lambda,\mu}=0$ unless $2k\in\N$, and for $2k\in\N$:
\begin{enumerate}[label=(\roman*)]
\item if $t_i\notin\N$ for some $i$, then $\dim\sT^n_{\lambda,\mu}=\binom{2k+n-2}{n-2}$;
\item if $t\in\N^n$, then
$$
\dim\sT^n_{\lambda,\mu}=\binom{2k+n-2}{n-2}+\max\big(0,\ \delta_{2k-1}(2t)-\delta_{2k}(2t)\big).
$$
\end{enumerate}
\end{theorem}

\begin{proof}
\emph{Step 1: reduction.} The field $X_1=\partial_x$ commutes with $D$, so it acts only on the coefficients: $X_1\cdot a_{m,\alpha}=m\,a_{m-1,\alpha}$ and $X_1\cdot b_{m,\alpha}=m\,b_{m-1,\alpha}$. An invariant operator therefore has coefficients independent of $x$, that is, it is a combination of the $a_{0,\alpha}=D^\alpha$ and $b_{0,\alpha}=\theta D^\alpha$. Moreover $X_x\cdot D^\alpha=(k-\frac{|\alpha|}2)D^\alpha$ and $X_x\cdot\theta D^\beta=(k+\frac12-\frac{|\beta|}2)\theta D^\beta$. By \eqref{gen}, an invariant $v$ is killed by $X_1$ and $X_x$. So $\sT^n_{\lambda,\mu}=0$ unless $2k\in\N$, and then
$$
v=\sum_{|\alpha|=2k}u_\alpha D^\alpha+\sum_{|\beta|=2k+1}w_\beta\,\theta D^\beta .
$$

\smallskip
\emph{Step 2: the equation $X_\theta\cdot v=0$.} The signs coming from the parity of the arguments are absorbed by the identity $\theta\,P(G)=\theta\,G$, where $P$ is the parity involution. With this, the computation of $X_\theta\cdot v=\frac12D\circ v-(-1)^{|v|}v\circ\frac12 D$ gives that $X_\theta\cdot v=0$ if and only if
\begin{equation}\label{E1}
w_\beta=-2\sum_{i\in\mathcal E(\beta)}u_{\beta-\varepsilon_i}\qquad(|\beta|=2k+1),
\end{equation}
the coefficients of $\theta D^\gamma$, $|\gamma|=2k+2$, vanishing then automatically. So $u$ is free and determines $w$. Note that this step does not involve $k$.

\smallskip
\emph{Step 3: the equation $X_{x\theta}\cdot v=0$.} From $D^2=\partial_x$ and $D(x)=\theta$ one gets, for every $H$ and $m\ge0$,
$$
D^m(xH)=xD^mH+\lfloor m/2\rfloor D^{m-2}H+\begin{cases}\theta\,D^{m-1}H,& m\ \text{odd},\\ 0,& m\ \text{even}.\end{cases}
$$
Using this (together with $\mathfrak L^{\lambda_i}_{X_{x\theta}}=\frac12xD+\lambda_i\theta$ and the same sign rules as in Step~2), a second computation shows that for $v$ as in Step~1
$$
X_{x\theta}\cdot v=x\,(X_\theta\cdot v)+Z,
$$
where $Z$ has constant coefficients, $Z=\sum_{|\beta|=2k-1}\zeta_\beta D^\beta+\sum_{|\gamma|=2k}\zeta'_\gamma\theta D^\gamma$, and
\begin{equation}\label{za}
\zeta_\beta=\sum_{i=1}^n(-1)^{s_i(\beta)}p_i(\beta+\varepsilon_i)\,u_{\beta+\varepsilon_i}.
\end{equation}
Assume $X_\theta\cdot v=0$, so that $X_{x\theta}\cdot v=Z$. By \eqref{gen},
$$
X_\theta\cdot Z=X_\theta X_{x\theta}\cdot v=\tfrac12X_x\cdot v-X_{x\theta}X_\theta\cdot v=0 .
$$
Since $Z$ has the shape of Step~1 with $k-\frac12$ in place of $k$, and Step~2 does not involve $k$, we get $\zeta'_\gamma=-2\sum_{i\in\mathcal E(\gamma)}\zeta_{\gamma-\varepsilon_i}$. Hence $Z=0$ if and only if $\zeta=0$. Together with Step~2:
$$
v\in\sT^n_{\lambda,\mu}\iff w\ \text{is given by \eqref{E1} and}\ \zeta_\beta=0\ \ (|\beta|=2k-1).
$$

\smallskip
\emph{Step 4: an odd operator.} Let $\Lambda=\R[z_1,\dots,z_n]\otimes\Lambda(\xi_1,\dots,\xi_n)$, with the weight $|a|+\frac12|\eta|$ on $z^a\xi^\eta$, and
$$
B=\sum_{i=1}^n\big((z_i\partial_{z_i}+2\lambda_i)\partial_{\xi_i}-\xi_i\partial_{z_i}\big).
$$
For $\alpha\in\N^n$ write $\alpha=2R+\eta$ with $\eta\in\{0,1\}^n$, put $f(\alpha)=\sum_j\alpha_j(\alpha_j-1)/2$ (so that $f(\beta+\varepsilon_i)-f(\beta)=\beta_i$) and
$$
\varphi_u=\sum_{|\alpha|=2k}(-1)^{f(\alpha)}u_\alpha\,z^R\xi^\eta ,
$$
an element of weight $k$; the map $u\mapsto\varphi_u$ is bijective onto the weight $k$ part of $\Lambda$. For $\beta=2R'+\eta'$ with $|\beta|=2k-1$, the coefficient of $z^{R'}\xi^{\eta'}$ in $B\varphi_u$ is
$$
2\,(-1)^{|\beta|+f(\beta)}\,\zeta_\beta .
$$
Indeed, $B$ sends $z^R\xi^\eta$ to the monomials $z^R\xi^{\eta-e_i}$ (when $\eta_i=1$, with coefficient $(-1)^{\eta_1+\dots+\eta_{i-1}}(R_i+2\lambda_i)$) and $z^{R-e_i}\xi^{\eta+e_i}$ (when $\eta_i=0$, $R_i\ge1$, with coefficient $-(-1)^{\eta_1+\dots+\eta_{i-1}}R_i$); in both cases the target is indexed by $\beta=\alpha-\varepsilon_i$, the sign $(-1)^{\eta_1+\dots+\eta_{i-1}}$ equals $(-1)^{|\beta|+s_i(\beta)}$, and $R_i+2\lambda_i=2p_i(\alpha)$ if $\alpha_i$ is odd, $R_i=2p_i(\alpha)$ if $\alpha_i$ is even. So the coefficient of $u_{\beta+\varepsilon_i}$ before the twist by $(-1)^f$ is $2(-1)^{|\beta|+s_i(\beta)+\beta_i}p_i(\beta+\varepsilon_i)$ (the factor $(-1)^{\beta_i}$ distinguishes the two cases), and $(-1)^{\beta_i}=(-1)^{f(\beta+\varepsilon_i)-f(\beta)}$ is exactly absorbed by the twist. Therefore
\begin{equation}\label{iso}
\sT^n_{\lambda,\mu}\ \simeq\ \{\varphi\in\Lambda\ \text{of weight }k:\ B\varphi=0\}.
\end{equation}

\smallskip
\emph{Step 5: chains.} $\Lambda=\bigotimes_iM_i$ with $M_i=\R[z_i]\otimes\Lambda(\xi_i)$, and $B$ acts by the graded Leibniz rule. On $M_i$ take the basis $e_a=z_i^a$, $e_{a+1/2}=z_i^a\xi_i$ ($a\in\N$); then
$$
B e_{a+1/2}=(a-t_i)\,e_a,\qquad Be_a=-a\,e_{a-1/2}.
$$
So $B$ lowers the weight by $\frac12$ and the only vanishing coefficients are $Be_0=0$ and, if $t_i\in\N$, $Be_{t_i+1/2}=0$. If $t_i\notin\N$, $M_i$ is one infinite chain $I^{(0)}$ starting at weight $0$. If $t_i\in\N$, the chain is cut and $M_i=J_i\oplus I_i$ as a module for $B$, with $J_i=\mathrm{Span}\{e_\nu:\nu\le t_i\}$ of dimension $2t_i+1$ and $I_i=\mathrm{Span}\{e_\nu:\nu\ge t_i+\frac12\}$ an infinite chain starting at weight $t_i+\frac12$. Put $A=\sum_i(z_i\partial_{\xi_i}-\xi_i)$, which raises the weight by $\frac12$; on $M_i$ one finds that the anticommutator of $A$ and $B$ is $-(2\nu-t_i)$ on the weight $\nu$ part, so $A$ and $B$ make $\Lambda$ an $\osp$-module. Then $I_i$ is an $\osp$-submodule (whereas $J_i$ is not: $A\,e_{t_i}=-e_{t_i+1/2}\in I_i$), and $M_i/I_i$, identified with $J_i$ as a $B$-module, is the irreducible module of dimension $2t_i+1$.

Two facts are needed. First, if $V$ is a graded space with an odd operator $B_V$ lowering the weight by $\frac12$, with weights bounded below and finite-dimensional weight spaces, and $I^{(b)}$ is an infinite chain starting at weight $b$, then the weight $\nu$ part of $\ker B$ on $V\otimes I^{(b)}$ has dimension $\dim V_{\nu-b}$. Indeed $\sum_jv_j\otimes e_{b+j/2}$ is killed by $B$ if and only if $B_Vv_j=\pm c_{j+1}v_{j+1}$ for all $j\ge0$ with $c_{j+1}\ne0$, so $v_0$ is free and determines the other $v_j$, and the sum is finite because the weights of $V$ are bounded below. Second, finite-dimensional $\osp$-modules are completely reducible, and on an irreducible one the kernel of the odd lowering operator $B$ is the lowest weight line.

\smallskip
\emph{Step 6: counting.} Let $q$ record the weight, $q^{2\nu}$ for weight $\nu$, so that $M_i$ has character $1/(1-q)$. Put $r_i=2t_i+1$ when $t_i\in\N$.

If $t_{i_0}\notin\N$, then $M_{i_0}=I^{(0)}$ and the first fact gives $\ker B\simeq\bigotimes_{j\ne i_0}M_j$, whose weight $k$ part has dimension $[q^{2k}](1-q)^{-(n-1)}=\binom{2k+n-2}{n-2}$. This is (i).

Let $t\in\N^n$. As in the classical case, for $S\subset\{1,\dots,n\}$ let
$$
\Lambda_S=\bigotimes_{i\in S}I_i\otimes\bigotimes_{i\notin S}J_i ,
$$
so that $\Lambda=\bigoplus_S\Lambda_S$ as graded spaces and $B$ preserves each $\Lambda_S$. For $S\ne\emptyset$ pick $i_0\in S$; the first fact, with $b=t_{i_0}+\frac12$, shows that $\ker B$ on $\Lambda_S$ has character
$$
\frac{\prod_{i\notin S}(1-q^{r_i})\prod_{i\in S}q^{r_i}}{(1-q)^{n-1}}.
$$
Since $\sum_S\prod_{i\notin S}(1-q^{r_i})\prod_{i\in S}q^{r_i}=1$, the sum over $S\ne\emptyset$ is $\big(1-\prod_i(1-q^{r_i})\big)/(1-q)^{n-1}$. Now $\prod_i(1-q^{r_i})/(1-q)^n=\sum_m\delta_m(2t)q^m$, so the weight $k$ part of this sum has dimension
$$
\binom{2k+n-2}{n-2}-\big(\delta_{2k}(2t)-\delta_{2k-1}(2t)\big).
$$
For $S=\emptyset$: the sum $\Lambda'=\bigoplus_{S\ne\emptyset}\Lambda_S=\sum_iI_i\otimes\bigotimes_{j\ne i}M_j$ is an $\osp$-submodule of $\Lambda$, because each $I_i$ is one, and the projection of $\Lambda$ onto $\Lambda_\emptyset$ along $\Lambda'$ commutes with $B$. Hence, as a $B$-module, $\Lambda_\emptyset\simeq\Lambda/\Lambda'\simeq\bigotimes_i(M_i/I_i)$, and the latter is an $\osp$-module, a tensor product of finite-dimensional irreducible modules, with character $\chi=\sum_m\delta_m(2t)q^m$, symmetric about $m=\tau=|t|$. It is completely reducible, and the irreducible component with lowest grade $m_0$ has character $q^{m_0}+\dots+q^{2\tau-m_0}$. If $c_{m_0}$ components have lowest grade $m_0$, then $(1-q)\chi=\sum c_{m_0}(q^{m_0}-q^{2\tau-m_0+1})$, so
$$
c_m=\delta_m(2t)-\delta_{m-1}(2t)\quad(m\le\tau),\qquad c_m=0\quad(m>\tau).
$$
By the second fact, $\ker B$ on $\Lambda_\emptyset$ has dimension $c_{2k}$ in weight $k$. Adding the two contributions: for $2k\le\tau$ the total is $\binom{2k+n-2}{n-2}$; for $2k>\tau$ it is $\binom{2k+n-2}{n-2}+\delta_{2k-1}(2t)-\delta_{2k}(2t)$, and this difference equals $c_{2\tau-2k+1}\ge0$ by symmetry. For $2k\le\tau$ we have $\delta_{2k-1}(2t)-\delta_{2k}(2t)=-c_{2k}\le0$. So both cases are given by the formula in (ii).
\end{proof}

\begin{remark}
(1) For $n=1$ the equations $\zeta_\beta=0$ can be written down directly. If $k=j+\frac12$, the only equation is $(j+2\lambda)u_{2j+1}=0$; if $k\in\N$, $k\ge1$, it is $u_{2k}=0$; if $k=0$ there is none. Hence $\dim\sT^1_{\lambda,\mu}=1$ for $k=0$ (the identity) and for $k=j+\frac12$ with $j=-2\lambda\in\N$, with generator $D^{2j+1}-2\,\theta D^{2j+2}$ by \eqref{E1}, and $\dim\sT^1_{\lambda,\mu}=0$ otherwise.\\
(2) Unlike the classical case, no Lefschetz argument is needed: the symmetry and unimodality of $\delta_m(2t)$ come from the representation theory of $\osp$.
\end{remark}

\begin{example}
Take $n=2$, $\lambda=(-\frac12,-\frac12)$, so $t=(1,1)$ and $\delta_m(2t)=(1,2,3,2,1)$. For $k=\frac32$ the equations \eqref{za} with $\beta=(2,0),(1,1),(0,2)$ read
$$
-\tfrac12u_{21}=0,\qquad \tfrac12u_{21}-\tfrac12u_{12}=0,\qquad -\tfrac12u_{12}=0 .
$$
So $u_{30}$ and $u_{03}$ are free and $\dim\sT^2_{\lambda,\mu}=2$, which agrees with $\binom30+\max(0,\delta_2(2t)-\delta_3(2t))=1+(3-2)$. For $k=2$ the equations for $|\beta|=3$ are $-u_{40}-\frac12u_{31}=0$, $u_{22}=0$ (twice) and $\frac12u_{13}-u_{04}=0$, so $u_{31}$ and $u_{13}$ are free and the dimension is again $2=1+(\delta_3(2t)-\delta_4(2t))$. For $k=1$ and any $\lambda$ the equations are $-\frac12u_{20}+\lambda_2u_{11}=0$ and $-\lambda_1u_{11}-\frac12u_{02}=0$, so $u_{11}$ is free and the dimension is $1$.
\end{example}

\end{document}